\documentclass[11pt,a4paper,oneside,leqno]{article}

\usepackage{tabularx}
\usepackage{enumitem}
\usepackage{amsmath,amssymb,amsthm,amsrefs,stmaryrd}
\usepackage[hidelinks]{hyperref}

\hypersetup{
    pdfauthor={Tom\'a\v{s} Crh\'ak},
    pdftitle={A note on sigma algebras on convex spaces},
    pdfkeywords={convex space, sigma algebra},
}

\setenumerate[1]{label=(\arabic*)}

\newtheorem{lemma}{Lemma}
\newcommand{\Set}{\textbf{Set}}
\newcommand{\Meas}{\textbf{Meas}}
\newcommand{\Cvx}{\textbf{Cvx}}
\newcommand{\eval}[1]{ev_{#1}}
\newcommand{\countables}[1]{[#1]^{\le\omega}}
\renewcommand{\complement}[1]{#1^{\sf{c}}}
\newcommand{\AND}{\:\&\:}
\newcommand{\I}{\mathbb{I}}
\newcommand{\A}{\mathcal{A}}
\newcommand{\E}{\mathcal{E}}
\newcommand{\C}{\mathcal{C}}
\renewcommand{\P}{\mathcal{P}}
\renewcommand{\S}{\mathcal{S}}
\newcommand{\sig}{$\sigma$}
\newcommand{\preimage}[1]{#1^{\shortleftarrow}}
\newcommand{\image}[1]{#1^{\shortrightarrow}}
\newcommand{\sbool}{\Sigma_{bool}}
\newcommand{\seval}{\Sigma_{ev}}
\newcommand{\restrict}[1]{{\big|}_{#1}}
\newcommand{\set}[1]{\left\{ #1 \right\}}
\newcommand{\ks}{\textsc{ks}}

\begin{document}

\title{A note on $\sigma$-algebras on sets of affine and measurable maps to the unit interval}
\author{Tom\'a\v{s} Crh\'ak}
\date{7 Mar 2018}
\maketitle

\begin{abstract}
In \textit{The factorization of the Giry monad} Kirk Sturtz
considers two $\sigma$-algebras on convex spaces of functions to the unit
interval. One of them is generated by the Boolean subobjects and the other
is the $\sigma$-algebra induced by the evaluation maps. The author asserts
that, under the assumptions given in the paper, the two $\sigma$-algebras
coincide. We give examples contradicting this statement.
\end{abstract}


\section{Introduction}

The notation and terminology of this paper is mostly based on
that of~\cite{sturtz2018factorization-v2}; $\Cvx$ is the category of convex
spaces with affine maps, $\Meas$ is the category of measurable spaces with
measurable maps and $\Set$ denotes the category of sets.

For a convex space $A$, the binary operations defining its convex structure
will be denoted by $+_r$, $r\in[0,1]$.  A subset $E\subseteq A$ is said to
be \textit{convex} if it is closed under all $+_r$ and it is said to be
a \textit{Boolean subobject} of $A$ if both $E$ and its complement
$\complement{E}$ are convex. It is immediate that the empty set, $A$
and also the inverse image of a convex set under an affine map are all
convex, and similarly for Boolean subobjects.

Depending on the context, the unit interval $\I$ will be viewed either
as an object of $\Meas$ or $\Cvx$. When considered as a measurable space,
the Borel \sig-algebra will be assumed. As a convex space, it is endowed
with the natural convex structure given by
$$
    x +_r y = (1-r)x + ry
    \qquad
    \forall r\in[0,1]
    \text{.}
$$
The Boolean subobjects of $\I$ are exactly the empty set and the convex sets
containing $0$ or $1$, i.e., the intervals
$[0, x)$, $[0, x]$, $(x, 1]$ and $[x, 1]$ for all $x\in\I$.


\section{The measurable structure of a convex space}

For a convex space $A$, the \textit{Boolean \sig-algebra} on $A$ is by
definition the smallest \sig-algebra containing all Boolean subobjects of $A$,
and is denoted by $\sbool^A$. In particular, the Boolean \sig-algebra
on $\I$ coincides with its Borel \sig-algebra. Note that for arbitrary
convex spaces $A$ and $B$, every affine map $A\to B$ is measurable
w.r.t.\ the Boolean \sig-algebras on $A$ and $B$.

For an arbitrary set $A$, $\Set(A,\I)$ may be turned into a convex space
by defining $+_r$ pointwise, thus for $f, g\in\Set(A,\I)$ and $a\in A$
we have
$$
    (f +_r g)(a) = f(a) +_r g(a)
    \qquad
    \forall r\in[0,1]
    \text{.}
$$
We are concerned with convex subsets $\A\subseteq\Set(A,\I)$.
In addition to the Boolean \sig-algebra, $\A$ may be endowed with the
\textit{evaluation \sig-algebra}, which is induced by the evaluation maps
$\eval{a}:\A\to\I$, $a\in A$, where $\eval{a}$ is defined by
$\eval{a}(f) = f(a)$. The evaluation \sig-algebra will be denoted
by $\seval^\A$. Since the evaluations are affine, it follows that
$$
    \seval^\A \subseteq \sbool^\A
    \text{.}
$$


\section{Overview of refuted statements}

Let $A$ be a convex space, and let $\Sigma A = \left( A, \sbool^A \right)$.
In~\cite{sturtz2018factorization-v2} the author
deals with the Boolean and evaluation \sig-algebras on $\Cvx(A,\I)$
and $\Meas(\Sigma A,\I)$, and asserts:

\begin{itemize}
\item Lemma 5.1.
    \textit{For any convex space $A$}
    \begin{equation*}\label{KS-5.1a}\tag{\ks-5.1a}
        \seval^{\Cvx(A,\I)} = \sbool^{\Cvx(A,\I)}
    \end{equation*}
    \begin{equation*}\label{KS-5.1b}\tag{\ks-5.1b}
        \seval^{\Meas(\Sigma A,\I)} = \sbool^{\Meas(\Sigma A,\I)}
    \end{equation*}
\item Corollary 5.2. \textit{For any affine map $P:\Meas(\Sigma A,\I)\to\I$}
    \begin{equation*}\label{KS-5.2}\tag{\ks-5.2}
        \textit{P is measurable w.r.t.\ $\seval^{\Meas(\Sigma A,\I)}$}
    \end{equation*}
\item Lemma 5.3.
    \textit{For any measurable space $X$}
    \begin{equation*}\label{KS-5.3}\tag{\ks-5.3}
        \seval^{\Meas(X,\I)} = \sbool^{\Meas(X,\I)}
    \end{equation*}
\end{itemize}

We will show that taking the unit interval for $A$ contradicts
(\ref{KS-5.2}), and hence (\ref{KS-5.1b}), and thus, for $X = \Sigma A$,
disproves (\ref{KS-5.3}) as well; it is, however, in accordance with
(\ref{KS-5.1a}). Next we show that the free convex space over an uncountable
set contradicts (\ref{KS-5.1a}) and (\ref{KS-5.1b}).


\section{Main results}

Throughout this section, $A$ is a convex space, and $\A$ is a subset of
$\Set(A,\I)$. For a function $f$ and a set $u$, $f\restrict{u}$,
$\image{f}(u)$ and $\preimage{f}(u)$ denote the restriction, image
and inverse image, respectively.

Let us define $\E\subseteq\P(\A)$ by
$$
    \E = \left\{
            E\subseteq\A:
                \big(\exists u\in\countables{A}\big)
                \big(\forall f\in E\big)
                \big(\forall g\in\A\big)
                \big(f\restrict{u} = g\restrict{u} \implies g\in E\big)
         \right\},
$$
where $\countables{A}$ denotes the set of countable subsets of $A$.
Then $\E$ has the following properties:

\begin{itshape}
\begin{enumerate}
\item $\emptyset\in\E$, $\A\in\E$;
\item $\E$ is closed under countable unions and intersections;
\item for $K\subseteq\I$ and $a\in A$, $\preimage{\eval{a}}(K)\in\E$.
\end{enumerate}
\end{itshape}

\begin{proof}
(1) is obvious.
Regarding (2), let $\C\subseteq\E$ be countable. For each $E\in\C$ fix $u_E$
such that
$$
    \big(\forall f\in E\big)
    \big(\forall g\in\A\big)
    \big(f\restrict{u_E} = g\restrict{u_E} \implies g\in E\big)
$$
and employ $u = \bigcup\{u_E: E\in\C\}$, which is countable,
to show that both $\bigcup\C$ and $\bigcap\C$ belong to $\E$.
The proof of (3) is obvious on setting $u = \{a\}$.
\end{proof}

Now it is immediate that
$\S = \left\{ E\subseteq\A: E\in\E \AND \complement{E}\in\E \right\}$
is a \sig-algebra and from
$\complement{(\preimage{\eval{a}}(K))} = \preimage{\eval{a}}(\complement K)$
it follows that $\preimage{\eval{a}}(K)\in\S$.
We conclude that $\seval^\A\subseteq\S$, which brings us to
\begin{lemma}\label{S}
For every convex space $A$, and every $\A\subseteq\Set(A,\I)$:
$$
    \seval^\A\subseteq\E
    \text{.}
$$
\end{lemma}


We have already mentioned that $\Sigma\I = \I$. From the following two
lemmas it follows that $A = \I$ provides a counterexample to (\ref{KS-5.1b}),
though it is in accordance with (\ref{KS-5.1a}).

\begin{lemma}\label{affine not measurable}
Let $P:\Meas(\I,\I)\to\I$ map every function of $\Meas(\I,\I)$ to its
Lebesgue integral. Then $P$ is affine but not measurable with respect to
the evaluation \sig-algebra on $\Meas(\I,\I)$.

\begin{proof}
We apply Lemma~\ref{S} for $A = \I$ and $\A = \Meas(\I,\I)$ to show that
$P^{-1}(0)\notin\seval^{\Meas(\I,\I)}$, from which the conclusion
of the lemma follows. For any countable subset $u\subseteq\I$,
consider the zero constant function $\overline{0}\in P^{-1}(0)$,
and the characteristic function
$\chi_{\complement{u}}\in\Meas(\I,\I)$. Then
$\overline{0}\restrict{u} = \chi_{\complement{u}}\restrict{u}$, but
$P(\chi_{\complement{u}}) = 1$, thus $\chi_{\complement{u}}\notin P^{-1}(0)$.
Therefore, $P^{-1}(0)\notin\E$.
\end{proof}
\end{lemma}


\begin{lemma}
The Boolean and evaluation \sig-algebras coincide for $\Cvx(\I,\I)$.

\begin{proof}
The affine map $\Phi:\Cvx(\I,\I)\to\I\times\I$ given by
$$
    \Phi(f) = (f(0), f(1))
$$
is an isomorphism of convex spaces. Under this isomorphism, the evaluations
$\eval{0}$ and $\eval{1}$ correspond to the projections of the product,
and it is straightforward that $\seval^{\I\times\I}$ coincides with the usual
Borel \sig-algebra on $\I\times\I$. Thus, given a Boolean subobject $E$
of $\I\times\I$, our business is to show that $E$ is a Borel set.

Recall that for every pair of elements $p, q\in\I\times\I$ there is
a unique affine map $\pi_{p,q}:\I\to\I\times\I$ such that
$\pi_{p,q}(0) = p$ and $\pi_{p,q}(1) = q$, which is given explicitly by
$$
    \pi_{p,q}(r) = p +_r q
    \text{.}
$$

From the Closed map lemma it follows that $\pi_{p,q}$ is a closed map, and from
this it is clear that the image of every convex subset of $\I$ under $\pi_{p,q}$
is a Borel subset of $\I\times\I$.

Let us now consider the set $\preimage{\pi_{p,q}}(E)$. As it is a Boolean
subobject of $\I$, its image under $\pi_{p,q}$ is a Borel set; we will
denote this image by $E_{p,q}$. In the sequel, we will employ the following
fact:
$$
    E_{p,q} \subseteq E
    \quad\AND\quad
    \image{\pi_{p,q}}(\I)\setminus E_{p,q} \subseteq \complement{E}
    \text{.}
$$

If $p\in E$, then $0\in\preimage{\pi_{p,q}}(E)$;
the supremum of $\preimage{\pi_{p,q}}(E)$ will be called
\textit{the dividing point determined by $\pi_{p,q}$}.
The dividing point determined by $\pi_{p,q}$, denoted for the moment by $u$,
is characterized by the following property:
$$
    \image{\pi_{p,q}}([0,u)) \subseteq E
    \quad\AND\quad
    \image{\pi_{p,q}}((u,1]) \subseteq \complement{E}
    \text{.}
$$

Now that we are through the preliminaries, we will distinguish four cases,
depending on how many corners of $\I\times\I$ belong to $E$ and $\complement{E}$.

\begin{enumerate}
\item If $E$ does not contain any corner, then all of them are contained
    in $\complement{E}$, and we have $\complement{E} = \I\times\I$, so that
    $E = \emptyset$, hence $E$ is a Borel set.

\item Suppose $E$ contains one corner and $\complement{E}$ contains the
    remaining three corners. We proceed for $(0,0)\in E$, the other cases
    are analogous. Let $u$ and $v$ be the dividing points determined
    by $\pi_{(0,0), (1,0)}$ and $\pi_{(0,0), (1, 0)}$, respectively.
    Then it is readily checked that
    $$
        E = \{ (x,y)\in\I\times\I: uy + vx < uv \}
            \cup
            E_{(0,v),(u,0)}
            \text{,}
    $$
    thus $E$ is a Borel set.

\item Suppose $E$ contains two \textit{neighboring} corners and $\complement{E}$
    contains the other two corners. We proceed for $(0,0), (1,0)\in E$,
    the other cases are analogous. Let $u$ and $v$ be the dividing points
    determined by $\pi_{(0,0), (0,1)}$ and $\pi_{(1,0), (1,1)}$, respectively.
    Then we have
    $$
        E = \{ (x,y)\in\I\times\I: y < (1-x)u + xv \}
            \cup
            E_{(0,u),(1,v)}
            \text{,}
    $$
    thus $E$ is a Borel set.

\item Suppose $E$ contains two \textit{opposite} corners. Then it must also
    contain at least one of the remaining corners, because otherwise the point
    $\left(\frac{1}{2}, \frac{1}{2}\right)$ would belong to both $E$
    and $\complement{E}$, a contradiction. But then we deal with a complementary
    case to (1) or (2).

\end{enumerate}
\end{proof}
\end{lemma}


The next result builds on a free convex space $A$ over an (arbitrary)
uncountable set $M$. Recall that $A$ is the set of maps $\alpha\in\Set(M,\I)$
with $\alpha^{-1}(0)$ cofinite in $M$ and $\sum_{m\in M}\alpha(m) = 1$.
$A$ is endowed with the pointwise convex structure.
As a consequence of the forthcoming lemma we have

\begin{enumerate}
\item $\seval^{\Cvx(A,\I)} \subsetneq \sbool^{\Cvx(A,\I)}$,
    which contradicts (\ref{KS-5.1a});

\item $\seval^{\Meas(\Sigma A,\I)} \subsetneq \sbool^{\Meas(\Sigma A,\I)}$,
    which contradicts (\ref{KS-5.1b}).
\end{enumerate}

\begin{lemma}\label{ev not bool}
Let $\overline{0}$ be the zero constant function $A\to\I$. For every convex
subset $\A\subseteq\Set(A,\I)$ that includes $\Cvx(A,\I)$ we have
$$
    \set{\overline{0}} \in \sbool^\A \setminus \seval^\A
    \text{.}
$$

\begin{proof}
Convexity of $\set{\overline{0}}$ and its complement is facile, we
focus on $\set{\overline{0}}\notin\seval^\A$. For a countable subset $u\subseteq A$
let $M_0 = \bigcap\set{\alpha^{-1}(0): \alpha\in u}$, and define
$g:A\to\I$ by
$$
    g(\alpha) = \sum_{m\in M_0}\alpha(m).
$$
Then $g$ is affine, thus $g\in\A$. For $\alpha\in u$ we have $g(\alpha) = 0$,
so that $g\restrict{u} = \overline{0}\restrict{u}$. However, $g\ne \overline{0}$.
To see this, observe that $M_0$ is non-empty as it is a countable intersection
of cofinite subsets of $M$, which has been chosen to be uncountable.
For an arbitrary $m_0\in M_0$ we take the unique $\alpha_0\in A$
with $\alpha_0(m_0) = 1$, hence $g(\alpha_0) = 1$. This proves
$g\notin\set{\overline{0}}$ and, as $u$ was arbitrary, we conclude that
$\set{\overline{0}}\notin\E$. Now apply Lemma~\ref{S} to show
$\set{\overline{0}}\notin\seval^\A$.
\end{proof}
\end{lemma}


\end{document}